\documentclass[10pt]{article}
\makeatletter
\@addtoreset{equation}{section}
\makeatother

\input epsf
\usepackage{amssymb}
\usepackage{amsmath}
\usepackage{epsfig}
\usepackage{pstricks}

\def\tfrac#1#2{{\frac{\lower.6ex
\hbox{$\scriptstyle#1$}}%\over
{\raise.7ex
\hbox{$\scriptstyle#2$}}}}

\def\bigO{{\cal O}}

\def\ZZ{{\mathbb Z}}
\def\calC{{\cal C}}

\def\dsp#1{\displaystyle#1}

\newcommand{\binomial}[2]{\left( \begin{array}{c} {#1} \\ 
                        {#2} \end{array} \right)} % binomial coefficients

\def\CHFs#1#2#3{
{}_1F_1\!\left({#1};{#2};{#3}\!\right)
}

\def\CHF#1#2#3{
{}_1F_1\left(
\begin{array}{c}
\begin{array}{c}\hskip-10pt#1\end{array}\\
\begin{array}{c}\hskip-10pt #2\end{array}
\end{array}
\hskip-8pt;\,#3
\right)}

\def\protectbold#1{\protect{\boldmath{$#1$}}}

\begin{document}

\title{Computing the Kummer function $U(a,b,z)$\\ for
  small values of the arguments}

\author{Amparo Gil\\
Departamento de Matem\'atica Aplicada y CC. de la Computaci\'on.\\
ETSI Caminos. Universidad de Cantabria. 39005-Santander, Spain.\\
\and
Javier Segura\\
        Departamento de Matem\'aticas, Estad\'{\i}stica y 
        Computaci\'on,\\
        Univ. de Cantabria, 39005 Santander, Spain.\\   
\and
Nico M. Temme\\
IAA, 1391 VD 18, Abcoude, The Netherlands\footnote{Former address: Centrum Wiskunde \& Informatica (CWI), 
        Science Park 123, 1098 XG Amsterdam,  The Netherlands}.\\ \\
{ \small
  e-mail: {\tt
    amparo.gil@unican.es,
    javier.segura@unican.es, 
    nico.temme@cwi.nl}}
    }

\date{\today}

\maketitle
\begin{abstract}
We describe methods for computing the Kummer function $U(a,b,z)$ for small values of $z$, 
with special attention to small values of $b$. For these values of $b$ the connection formula 
that represents $U(a,b,z)$ as a linear combination of two ${}_1F_1$-functions needs a limiting 
procedure. We use the power series of the ${}_1F_1$-functions and consider the terms for which this
limiting procedure is needed. We give recursion relations for higher terms in the expansion, and we 
consider the derivative $U^\prime(a,b,z)$ as well.
We also discuss the performance for small $\vert z\vert$ of an asymptotic approximation
 of the Kummer function in terms of modified Bessel functions.
\end{abstract}

\vskip 0.8cm \noindent
{\small
\par\noindent
Keywords \& Phrases:
Kummer function; numerical computation.
}

\section{Introduction}\label{sec:intro}
As is well known, confluent hypergeometric functions appear in a vast number of applications in physics and engineering. In spite of
their importance, few algorithms are available for the computation of any of
the standard solutions of the Kummer equation in the case of real or complex parameters; see for example
\cite{Muller:2001:CCH,Nardin:1992:ACM} for the computation of $M(a,b,x)$.  
In this paper we describe a method for computing the Kummer function $U(a,b,z)$ for small values of $\vert z\vert$, $\vert a\vert$ and $\vert b\vert$
  by using a simple expansion of a function in terms of reciprocal gamma functions.
To compute the Kummer function $U(a,b,z)$ for small values of $\vert z\vert$, 
the connection formulas for this function in terms of the ${}_1F_1$-functions is the most obvious starting point. When also small values of the $b$-parameter are allowed,  special care is needed because certain terms in these formulas will become singular, although the limit for $b\to0$ will exist. In case $b=0$ the power series expansion of $U(a,b,z)$ will contain a logarithmic term. This phenomenon occurs as well for positive integer values of $b$, but when the case for $b=0$ is under control, we can use stable recursions \cite{Segura:2008:NSS} to handle $b=1,2,3,\ldots$. 
It should also be noted that, for numerical evaluations with fixed precision, the power series expansions can only be used for rather small positive  values of $z$, say, for $0< z\le 1$. This is due to  the behaviour of the ${}_1F_1$-functions, which become exponentially large as $z\to +\infty$, whereas the $U$-function is only of algebraic growth in that case.

 With the coefficients provided 
by the expansion in terms of reciprocal gamma functions,
we can safely compute the function $U(a,b,z)$ and its $z$-derivative in double precision for 
\begin{equation}\label{eq:kumint01}
0< \vert z\vert\le 1, \quad \vert a\vert \le \tfrac12, \quad \vert b\vert \le \tfrac12.
\end{equation}
This computational method complements others described in \cite{Temme:1983:NCC}, where the computation of the confluent
hypergeometric function $U(a,b,z)$ is considered in other (complementary) parameter domains.

In the paper we also consider an asymptotic approximation 
of the Kummer function  $U(a,b,z)$ in terms of modified $K$-Bessel functions
and discuss its numerical performance for small $\vert z\vert$. An analysis of the stability of the recurrence
relations satisfied by the coefficients of the expansion is provided as well as few tests illustrating its
behaviour.

\section{Power series expansions for small \protectbold{\vert z\vert}}\label{sec:smallz}
We use the representation (see \cite[Eqn.~13.2.E41]{Olde:2010:CHF})
\begin{equation}\label{eq:kumsz01}
U(a,b,z)=\frac{\Gamma(1-b)}{\Gamma(a-b+1)}\CHF{a}{b}{z}+
\frac{z^{1-b}\Gamma(b-1)}{\Gamma(a)}\CHF{a-b+1}{2-b}{z}
\end{equation}
and the form
\begin{equation}\label{eq:kumsz02}
U(a,b,z)=\frac{\pi}{\sin(\pi b)}\left(\frac{\CHFs{a}{b}{z}}{\Gamma(b)\Gamma(a-b+1)}-
z^{1-b}\frac{\CHFs{a-b+1}{2-b}{z}}{\Gamma(a)\Gamma(2-b)}\right).
\end{equation}
These connection formulas are not defined for integer values of $b$, although the limit exists for, say, $b\to 0$. 

We assume $z\ne0$ and $b\notin\ZZ$. We expand the ${}_1F_1$-functions and obtain
\begin{equation}\label{eq:kumsz03}
U(a,b,z)=\frac{\Gamma(1-b)}{\Gamma(a-b+1)}+\frac{\pi b\,z}{\sin(\pi b)\Gamma(a)\Gamma(a-b+1)}
\sum_{m=0}^\infty w_m \frac{z^m}{m!},
\end{equation}
where
\begin{equation}\label{eq:kumsz04}
\begin{array}{@{}r@{\;}c@{\;}l@{}}
w_m&=&\dsp{\frac{u_m}{v_m}},\\[8pt]
u_m&=&(A_m-B_m)/b,\\[8pt]
A_m&=&\Gamma(m+1)\Gamma(2-b+m)\Gamma(a+m+1),\\[8pt]
B_m&=& z^{-b}\,\Gamma(a-b+1+m)\Gamma(b+1+m)\Gamma(m+2),\\[8pt]
v_m&=&\,\Gamma(m+2)\Gamma(b+m+1)\Gamma(2-b+m).
\end{array}
\end{equation}

The $u_m$ are well defined in the limit $b\to0$, but for numerical computations we need a stable representation of $w_m$. We derive such a representation that we can use for values of $a$ and $b$ in the interval $[-\frac12,\frac12]$. For other $a$-values we need a few recursion steps when this parameter is not too large.  For recursion with respect to $b$, we refer to \S\ref{sec:posintb}.

We concentrate on  $w_0$; other values $w_m$ follow from stable recursions  when $b$ is small. Namely, it is easily verified that
\begin{equation}\label{eq:kumsz05}
\begin{array}{@{}r@{\;}c@{\;}l@{}}
u_{m+1}&=& (a_m A_m-b_m B_m)/b=a_mu_m+d_mB_m,\\[8pt]
 d_m&=&(a_m-b_m)/b,
\end{array}
\end{equation}
where 
\begin{equation}\label{eq:kumsz06}
\begin{array}{@{}r@{\;}c@{\;}l@{}}
a_m&=&c_m-(m^2+(a+2)m+a+1)b,\\[8pt]
b_m&=&c_m+(m+2)(a-b)b,\\[8pt]
c_m&=&(m+1)(m+2)(m+a+1),\\[8pt]
d_m&=&-(m^2+2m(a+1)+3a+1)+(m+2)b.
\end{array}
\end{equation}
Theorem 4.17 of \cite{Gil:2007:NSF} can be used to check that $\{u_m\}$ is a dominant solution of (\ref{eq:kumsz05}), 
which means that when we have a stable representation of $u_m$, we can compute $u_{m+1}$ in a stable way as $b\to0$.
The recursion for $v_m$ is straightforward, with starting value
\begin{equation}\label{eq:kumsz07}
v_0=\Gamma(b+1)\Gamma(2-b)=\frac{(1-b)\,\pi b}{\sin(\pi b)}.
\end{equation}

\subsection{Computing the derivative}\label{sec:deriv}
For the derivative of $U(a,b,z)$ we obtain from \eqref{eq:kumsz03}
\begin{equation}\label{eq:kumsz08}
U^\prime(a,b,z)=\frac{\pi b}{\sin(\pi b)\Gamma(a)\Gamma(a-b+1)}
\sum_{m=0}^\infty w_m^d \frac{z^m}{m!},
\end{equation}
where
\begin{equation}\label{eq:kumsz09}
w_m^d=(m+1) w_m+zw_m^\prime.
\end{equation}
For $w_m^\prime$ we need, using \eqref{eq:kumsz04}, 
\begin{equation}\label{eq:kumsz10}
u_m^\prime=\frac{1}{z}B_m,
\end{equation}
for which we do not need further analysis.

\subsection{Computing the first term  \protectbold{w_0}}\label{sec:w0}
Following the ideas described in \cite{Skorokhodov:2001:ARM}, we introduce the function
\begin{equation}\label{eq:kumw001}
G(a,b)=\frac{1}{b}\left(\frac{1}{\Gamma(a+1+b)}-\frac{1}{\Gamma(a+1)}\right).
\end{equation}
The reciprocal gamma function is an entire function, and we have the expansion
\begin{equation}\label{eq:kumw002}
\frac{1}{\Gamma(z)}=\sum_{k=1}^\infty c_k z^k, \quad \vert z\vert <\infty,
\end{equation}
where $c_1=1$, $c_2=\gamma$, and the rest of the coefficients satisfy (see \cite[\S5.7]{Askey:2010:GFN})
\begin{equation}\label{eq:kumw003}
(k-1)c_k=\gamma c_{k-1}-\zeta(2)c_{k-2}+\zeta(3)c_{k-3}-\ldots+(-1)^k\zeta(k-1)c_1,\quad k\ge3,
\end{equation}
where $\gamma$ is Euler's constant and  $\zeta(z)$ denotes the $\zeta$-function. See Table~\ref{table:coefck}; for 31D values see \cite{Wrench:1968:CTS}.

For $G(a,b)$ defined in \eqref{eq:kumw001} we have the representation
\begin{equation}\label{eq:kumw004}
G(a,b)=\sum_{k=2}^\infty c_k d_k,\quad d_k=\frac{1}{b}\left((a+b)^{k-1}-a^{k-1}\right),
\end{equation}
which we use for  $a$ and $b$ in the interval $[-\frac12,\frac12]$, or for complex values with modulus not exceeding $\frac12$. We have $d_2=1$,  $d_3=2a+b$, and the recurrence relation
\begin{equation}\label{eq:kumw005}
d_{k+2}=(2a+b)d_{k+1}-a(a+b)d_k,\quad k\ge1.
\end{equation}

\renewcommand{\arraystretch}{1.1}
\begin{table}%[H]
  \caption{Coefficients $c_k$ of the expansion in \eqref{eq:kumw002}. 
  \label{table:coefck}} 
  \begin{centering}
  {%\tt 
  {%\small
    \begin{tabular}{ll}
  & \\ \hline
\ $c_{ 1}$ = \ 1.00000000000000000000;        &  $c_{ 15}$ = -0.20563384169776e-6; \\
\  $c_{ 2}$ = \ 0.57721566490153286061;        &    $c_{ 16}$ =\  \ 0.611609510448e-8;\\
\  $c_{ 3}$ = -0.65587807152025388108;       &   $c_{ 17}$ =\  \ 0.500200764447e-8;\\
\  $c_{ 4}$ = -0.4200263503409523553e-1;    & $c_{ 18}$ = -0.118127457049e-8;\\
\  $c_{ 5}$ = \  0.16653861138229148950;         & $c_{ 19}$ =\  \ 0.10434267117e-9;\\
\  $c_{ 6}$ = -0.4219773455554433675e-1;    & $c_{ 20}$ =\  \ 0.778226344e-11;\\
\  $c_{ 7}$ = -0.962197152787697356e-2;      & $c_{ 21}$ = -0.369680562e-11;\\
\  $c_{ 8}$ = \ 0.721894324666309954e-2;       & $c_{ 22}$ =\  \ 0.51003703e-12;\\
\  $c_{ 9}$ = -0.116516759185906511e-2;      &  $c_{ 23}$ = -0.2058326e-13;\\
 $c_{ 10}$ = -0.21524167411495097e-3;      &    $c_{ 24}$ = -0.53481e-14;\\
 $c_{ 11}$ = \  0.12805028238811619e-3;       &  $c_{ 25}$ =\  \ 0.12268e-14;\\
 $c_{ 12}$ = -0.2013485478078824e-4;       & $c_{ 26}$ = -0.11813e-15;\\
 $c_{ 13}$ = -0.125049348214267e-5;         &   $c_{ 27}$ =\  \ 0.119e-17;\\
 $c_{ 14}$ = \  0.113302723198170e-5;          & $c_{ 28}$ =\  \ 0.141e-17;\\
\hline
    \end{tabular}} \\}
  \end{centering}
  \end{table}
\renewcommand{\arraystretch}{1.0}

For a few recursion steps we can use the relations
\begin{equation}\label{eq:kumw006}
\begin{array}{@{}r@{\;}c@{\;}l@{}}
(a+1)(a+b+1)G(a+1,b)&=&\dsp{(a+1)G(a,b)-\frac{1}{\Gamma(a+1)},}\\[8pt]
(a+2)(a+b+2)G(a+2,b)&=&(2a+b+3)G(a+1,b)-G(a,b).
\end{array}
\end{equation}
In this case, Theorems 4.6 and 4.7 of \cite{Gil:2007:NSF} are inconclusive with respect to the existence of minimal solutions for
these recurrences for $G(a,b)$.

For $w_0$ we have
\begin{equation}\label{eq:kumw007}
bw_0=\frac{\Gamma(a+1)}{\Gamma(b+1)}-z^{-b}\frac{\Gamma(a-b+1)}{\Gamma(2-b)}.
\end{equation}
We write
\begin{equation}\label{eq:kumw008}
bw_0=\frac{\Gamma(a+1)}{\Gamma(b+1)}-\frac{\Gamma(a-b+1)}{\Gamma(2-b)}-\left(z^{-b}-1\right)\frac{\Gamma(a-b+1)}{\Gamma(2-b)},
\end{equation}
and we express some of the gamma functions in terms of the function $G(a,b)$. We have
\begin{equation}\label{eq:kumw009}
\begin{array}{@{}r@{\;}c@{\;}l@{}}
\Gamma(a-b+1)&=&\dsp{\frac{\Gamma(a+1)}{1-b\Gamma(a+1)G(a,-b)},}\\[8pt]
\Gamma(1+b)&=&\dsp{\frac{1}{1+bG(0,b)},}\\[8pt]
\Gamma(2-b)&=&\dsp{\frac{1-b}{1-bG(0,-b)}.}
\end{array}
\end{equation}
and obtain
\begin{equation}\label{eq:kumw010}
\begin{array}{@{}r@{\;}c@{\;}l@{}}
w_0&=&\dsp{\frac{\Gamma(a+1)}{(b-1)\left(1-b\Gamma(a+1)G(a,-b)\right)}\ \times}\\[10pt]
&&\Bigl(\Gamma(a+1)(1-b)\bigl(1+bG(0,b)\bigr)G(a,-b)\ +\\[4pt]
&&\dsp{1+(b-1)G(0,b)-G(0,-b)+\frac{z^{-b}-1}{b}}\bigl(1-bG(0,-b)\bigr)\Bigr).
\end{array}
\end{equation}
When we use an expansion for
\begin{equation}\label{eq:kumw011}
\frac{z^{-b}-1}{b}=\ln z\frac{e^{-b\ln z}-1}{b\ln z},
\end{equation}
for small values of $\vert b\ln z\vert$, we have obtained a representation for $w_0$ in which we can allow small values of $b$.
To start the recursion for $u_m$ given in \eqref{eq:kumsz05} we use $u_0=v_0w_0$ with $v_0$ given in \eqref{eq:kumsz07}.

In \cite{Skorokhodov:2001:ARM} the function
\begin{equation}\label{eq:kumw012}
\Gamma_\epsilon(z)=\frac{1}{\epsilon}\left(\frac{\Gamma(z+\epsilon)}{\Gamma(z)}-1\right)
\end{equation}
is considered for small values of $\epsilon$. In our notation we have
\begin{equation}\label{eq:kumw013}
\Gamma_\epsilon(z)=-\Gamma(z+\epsilon)G(z-1,\epsilon).
\end{equation}
Skorokhodov introduced this function for the computation of the Gauss hypergeometric function ${}_2F_1(a,b;c;z)$, where for certain values of the parameters $a$, $b$ and $c$ numerical problems arise. Because Skorokhodov needs the evaluation of $\Gamma_\epsilon(z)$ for general complex values of $z$, more details and expansions are considered in  \cite{Skorokhodov:2001:ARM}. For other papers on  the computation of the ${}_2F_1$-function (with attention to the computation of functions like $\Gamma_\epsilon(z)$ and $G(a,b)$), we refer to \cite{Doornik:2014:NEG} and  \cite{Forrey:1997:CHF}

 \subsection{Positive integer values of \protectbold{b}}\label{sec:posintb}
 When we have control of the computations for small values of $b$, we can obtain the functions $U(a,b+k,z)$ without further extra analysis, although the connection formulas in \eqref{eq:kumsz01} and \eqref{eq:kumsz02} cannot be used for $b\in\ZZ$.
 
 We have the relations
 \begin{equation}\label{eq:kumw014}
\begin{array}{@{}r@{\;}c@{\;}l@{}}
U(a,b+1,z) &=& U(a,b,z)-U^\prime(a,b,z),\\[8pt]
zU^\prime(a,b+1,z) &=& bU^\prime(a,b,z)-aU(a,b,z),
\end{array}
\end{equation}
and we conclude that, when $b$ is small,  to obtain the values on the left sides no special limits are needed, once we have  $U(a,b,z)$ and its derivative. In addition, when $a$, $b$ and $z$ are positive, we have $U^\prime(a,b,z)<0$, and the recursion is perfectly stable. For more computational aspects of the recursion with respect to $b$, we refer to \cite{Temme:1983:NCC}.

\subsection{Computing \protectbold{G(a,b)} by quadrature}\label{sec:kumtest}
Because $1/\Gamma(z)$ is an entire function, we can write the function $G(a,b)$ defined in \eqref{eq:kumw001} in the form
\begin{equation}\label{eq:kumgab01}
G(a,b)=\frac{1}{2\pi i}\int_\calC \frac{1}{\Gamma(z+1)}\  \frac{dz}{(z-a)(z-a-b)},
\end{equation}
where $\calC$ is a closed contour that encircles  the points $z=a$ and $z=a+b$ clockwise, say a circle around the origin with radius $r>\max(\vert a\vert,\vert a+b\vert)$. Using such a circle, we can write
\begin{equation}\label{eq:kumgab02}
G(a,b)=\frac{1}{2\pi }\int_{-\pi}^\pi\frac{1}{\Gamma(z)}\  \frac{d\theta}{(z-a)(z-a-b)},\quad z=r e^{i\theta}.
\end{equation}
As explained in \cite[\S5.2.3]{Gil:2007:NSF} and \cite{Trefethen:2014:TEC}, the trapezoidal rule gives excellent results in a quadrature algorithm. For an efficient code, we can precompute the values of $1/\Gamma(z)$ in a set of equidistant nodes on the circle with radius $r$.

\subsection{A numerical test}\label{sec:w0q}
For a numerical test we have verified the relation
\begin{equation}\label{eq:kumt01}
U(a-1,b,z)=(a-b+z) U(a,b,z)-z U^\prime(a,b,z)
\end{equation}
for $a=0.2$ and a few values of $z$ and $b=10^{-2k}$, $k=1(1)5$. In Table~\ref{table:test} we give the relative error in this relation. We have used a Maple code with Digits=16 and the coefficients in Table~\ref{table:coefck}; we stop the summation of the series in \eqref{eq:kumsz03} and \eqref{eq:kumsz08} when the absolute values of the terms are less than $10^{-16}$. For $z= -0.5-0.1 i$ about 15 terms are needed, for $z=1+i$ about 20. Because of the strong relation between the terms in the series, the evaluation of the series is done together.

\renewcommand{\arraystretch}{1.1}
\begin{table}%[H]
  \caption{Relative errors in the relation in \eqref{eq:kumt01} for $a=0.2$ and $b=10^{-2k}$, $k=1(1)5$, and two complex values of $z$. 
  \label{table:test} } 
  \begin{centering}
  {%\tt 
  {%\small
    \begin{tabular}{lcc}
 &&   \\
  \hline
  $k$ &  $z= -0.5-0.1 i$ & $z=1+i$  \\
  \hline
1 & 0.15e-15  &  0.25e-15 \\
2 & 0.10e-14  & 0.76e-15 \\
3 & 0.11e-14  & 0.26e-14 \\
4 & 0.33e-15  & 0.23e-14 \\
5 & 0.10e-14  & 0.11e-14 \\
\hline
    \end{tabular}} \\}
  \end{centering}
  \end{table}
\renewcommand{\arraystretch}{1.0}

\section{Asymptotic representation in terms of modified Bessel functions}\label{sec:asrep}

Our starting point is the representations for $a$ large given in \cite[\S10.3]{Temme:2014:AMI} 

\begin{equation}\label{eq:asrep01}
\begin{array}{ll}
\dsp{U(a,b,z)= 2\left(\frac{z}{a}\right)^{\frac12(1-b)}\frac{e^{\frac12z}}{\Gamma(a)}\ \times}\\[8pt]
\quad\quad\dsp{\left(K_{b-1}\left(2\sqrt{az}\right)A_a(z,b)+\sqrt{\frac{z}{a}}K_{b}\left(2\sqrt{az}\right)B_a(z,b)\right),}
\end{array}
\end{equation}

and
\begin{equation}\label{eq:asrep02}
\begin{array}{ll}
\dsp{\frac{1}{\Gamma(b)}\CHF{a}{b}{z}=\left(\frac{z}{a}\right)^{\frac12(1-b)}\frac{\Gamma(1+a-b)e^{\frac12z}}{\Gamma(a)}\ \times}\\[8pt]
\quad\quad\dsp{\left(I_{b-1}\left(2\sqrt{az}\right)A_a(z,b)-\sqrt{\frac{z}{a}}I_{b}\left(2\sqrt{az}\right)B_a(z,b)\right),}
\end{array}
\end{equation}

where
\begin{equation}\label{eq:asrep03}
A_a(z,b)\sim \sum_{k=0}^\infty \frac{a_k(z)}{a^k}, \quad B_a(z,b)\sim \sum_{k=0}^\infty \frac{b_k(z)}{a^k},\quad a\to \infty,
\end{equation}
valid for bounded values of $\vert z\vert$ and $\vert b\vert$. The coefficients of these expansions can be computed by a simple 
scheme.

Let
\begin{equation}\label{eq:asrep04}
f(z,s)=e^{zg(s)}\left(\frac{s}{1-e^{-s}}\right)^c,\quad g(s)=\frac{1}{s}-\frac{1}{e^s-1}-\frac12.
\end{equation}
The function $f$ is analytic in the strip $\vert\Im s\vert<2\pi$ and it can be expanded for  $\vert s\vert<2\pi$ into
\begin{equation}\label{eq:kumUap05}
f(z,s)=\sum_{k=0}^{\infty} c_k(z) s^k.
\end{equation}
Then, the coefficients $a_k(z)$ and $b_k(z)$ of \eqref{eq:asrep03} are given by\footnote{Corrected form of \cite[Eq.~(10.3.41)]{Temme:2014:AMI}.}
\begin{equation}\label{eq:asrep06}
\begin{array}{@{}r@{\;}c@{\;}l@{}}
a_k(z) & = & \dsp{\sum_{m=0}^k \binomial{k}{m}(m+1-b)_{k-m}z^m c_{k+m}(z),}\\[8pt]
b_k(z) & = & \dsp{\sum_{m=0}^k \binomial{k}{m}(m+2-b)_{k-m}z^m c_{k+m+1}(z).}
\end{array}
\end{equation}

\subsection{Slater's results for large \protectbold{a}}\label{sec:Slater}
Slater's asymptotic expansions of the Kummer functions for large $a$ are given in \cite[\S4.6.1]{Slater:1960:CHF}, and are in terms of the large parameter 
$a$ written in the form
\begin{equation}\label{eq:Sla01}
a=\tfrac14u^2+\tfrac12b,
\end{equation}
where $u>0$ if $a$ and $b$ are real with $a>\frac12b$. For programming purposes we have simple recurrence relations for the coefficients in the asymptotic expansions.

Slater's results for large $a$ are written in the form
\begin{equation}\label{eq:Sla02}
\begin{array}{ll}
\dsp{e^{-\frac12z^2}z^b\CHF{a}{b}{z^2}\sim\Gamma(b)u^{1-b}2^{b-1}}\ \times \\[8pt]
\quad\quad\quad
\dsp{\left(zI_{b-1}(uz)\sum_{k=0}^{\infty}\frac{A_k(z)}{u^{2k}}+
\frac{z}{u}I_{b}(uz)\sum_{k=0}^{\infty}\frac{B_k(z)}{u^{2k}}\right)},
\end{array}
\end{equation}
and\footnote{Slater's result for the $U$-function contains an error; see \cite{Temme:2013:RSA} and \cite[\S10.3.6]{Temme:2014:AMI}.}
\begin{equation}\label{eq:Sla03}
\begin{array}{ll}
\dsp{e^{-\frac12z^2}z^bU\left(a,b,z^2\right)\sim\frac{2^{b}u^{1-b}}{\Gamma(1+a-b)}}\ \times \\[8pt]
\quad\quad\quad
\dsp{\left(zK_{b-1}(uz)\sum_{k=0}^{\infty}\frac{A_k(z)}{u^{2k}}-
\frac{z}{u}K_{b}(uz)\sum_{k=0}^{\infty}\frac{B_k(z)}{u^{2k}}\right)}.
\end{array}
\end{equation}
The coefficients are given by $A_0=1$ and
\begin{equation}\label{eq:Sla04}
\begin{array}{ll}
\dsp{B_k(z)=-\tfrac12A_k^\prime(z)+\int_0^z\left(\tfrac12t^2A_k(t)-\frac{b-\frac12}{t}A_k^\prime(t)\right)\,dt,} \\[8pt]
\dsp{A_{k+1}(z)=\frac{b-\frac12}{z}B_k-\tfrac12B_k^\prime(z)+\tfrac12\int^z t^2B_k(t)\,dt+K_k},
\end{array}
\end{equation}
and $K_k$ is chosen so that $A_{k+1}(z)\to0$ as $z\to0$. 
In fact,
\begin{equation}\label{eq:Sla05}
\begin{array}{ll}
\dsp{A_0(z)=1,\quad B_0(z)=\tfrac16z^3,}\\[8pt]
\dsp{A_1(z)=\tfrac16(b-2)z^2+\tfrac{1}{72}z^6,}\\[8pt]
\dsp{B_1(z)=-\tfrac13b(b-2)z-\tfrac{1}{15}z^5+\tfrac{1}{1296}z^{9},}\\[8pt]
\dsp{A_2(z)=-\tfrac{1}{120}(5b-12)(b+2)z^4+\tfrac{1}{6480}(5b-52)z^8+\tfrac{1}{31104}z^{12},}\\[8pt]
\dsp{B_2(z)=\tfrac{1}{90}(5b-12)(b+2)(b+1)z^3\ -}\\[8pt]
\quad\quad\quad\quad \dsp{\tfrac{1}{45360}(175b^2-350b-1896)z^7
-\tfrac{7}{12960}z^{11}+\tfrac{1}{933120}z^{15}.}
\end{array}
\end{equation}

\section{Convergent expansions based on the asymptotic representation}\label{sec:convrep}
We consider again the representation given in \eqref{eq:asrep01}, but instead of using asymptotic expansions of  $A_a(z,b) $ and $B_a(z,b) $ in negative powers of $a$,  we expand
\begin{equation}\label{eq:convrep01}
A_a(z,b) = \sum_{k=0}^\infty {\alpha_k(a,b)}{z^k}, \quad B_a(z,b)= \sum_{k=0}^\infty {\beta_k(a,b)}{z^k},
\end{equation}
where the coefficients follow from the recurrence equations ($k\ge 1$)
\begin{equation}\label{eq:convrep02}
\left\{
\begin{array}{ll}
\alpha_{k-1}=2b\alpha_{k}+4(k+1)(k+b)\alpha_{k+1}-4(2k+1)\beta_{k},\\ [8pt]
\beta_{k-1}=2b\beta_{k}+4(k+1)(k+2-b)\beta_{k+1}-8a(k+1)\alpha_{k+1}.
\end{array}
\right.
\end{equation}
These follow from the system of differential equations
\begin{equation}\label{eq:convrep03}
\left\{
\begin{array}{ll}
4zA^{\prime\prime}+4bA^{\prime}-8zB^{\prime}+(2b-z)A-4B=0,\\ [8pt]
4zB^{\prime\prime}+4(2-b)B^{\prime}-8aA^{\prime}+(2b-z)B=0,
\end{array}
\right.
\end{equation}
which can be obtained by applying Kummer's differential equation  
to the representations in \eqref{eq:asrep01} and \eqref{eq:asrep02}.

By solving \eqref{eq:asrep01} and \eqref{eq:asrep02} for the $A$- and $B$-functions, and replacing the 
$K$-Bessel functions by $I$-Bessel functions, we find that  $A_a(z,b)$ and $B_a(z,b)$ are entire functions of $z$. Hence, 
the series in \eqref{eq:convrep01} converge for all finite $z$.

First values are
\begin{equation}\label{eq:convrep04}
\begin{array}{ll}
\dsp{\alpha_0=\frac{a^{1-b}\Gamma(a)}{\Gamma(a+1-b)}, }& \dsp{\alpha_1=\frac{\alpha_0\left(b^2-b+2a\right)-2a}{2b(1-b)},}\\[8pt]
\dsp{\beta_0=\frac{a(\alpha_0-1)}{1-b}},  &\dsp{\beta_1=\frac{a\left(\alpha_0\left(4a-2b+b^2\right)-4a+b^2\right)}{2b(b-1)(b-2)}}.
\end{array}
\end{equation}
These values are well defined for $b=0,1,2$.

For large $a$ we can use  \cite[Eq.~(6.5.72)]{Temme:2014:AMI} for $\alpha_0$ and we see that  $\beta_0=1+\bigO(1/a)$. 
From a few Maple experiments, 
we see that a set of other coefficients in  the expansions in \eqref{eq:convrep01} are $\bigO(1/a)$, and it is 
likely that we should not use forward recursion, when $a$ is large.

We can change the mixed recursions in \eqref{eq:convrep02} into two recursion for $\alpha_k$ and $\beta_k$. We have ($k\ge3$)
\begin{equation}\label{eq:convrep05}
\begin{array}{ll}
p_{-3}\alpha_{k-3}+p_{-2}\alpha_{k-2}+p_{-1}\alpha_{k-1}+p_{0}\alpha_{k}+p_{1}\alpha_{k+1}=0,\\[8pt]
q_{-3}\beta_{k-3}+q_{-2}\beta_{k-2}+q_{-1}\beta_{k-1}+q_{0}\beta_{k}+q_{1}\beta_{k+1}=0,
\end{array}
\end{equation}
where
\begin{equation}\label{eq:convrep06}
\begin{array}{ll}
p_{-3} = -(2k-1)(2k+1),\\[8pt]
p_{-2} =8b(2k+1)(k-1),\\[8pt]
p_{-1} =-8+4b+8k^2+24k-64k^3+32k^4+12b^2-16kb+16bk^2-16k^2b^2+ 16kb^2,\\[8pt]
p_{0}= 16k(2k-3)(8ak^2-2bk^2-b^2-2a+b),\\[8pt]
p_{1} =16k(2k-1)(2k-3)(k+1)(k+b)(-k-1+b),
\end{array}
\end{equation}
and
\begin{equation}\label{eq:convrep07}
\begin{array}{ll}
q_{-3} =k, \\[8pt]
q_{-2} =-2b(2k-1),\\[8pt]
q_{-1} = 4(k-1)(-2k^2+b^2),\\[8pt]
q_{0}=  -8k(2k+1)(k-1)(4a-b),\\[8pt]
q_{1} = 16k(k-1)(k+1)(k+b)(k+2-b).
\end{array}
\end{equation}

In order to decide which of the directions of recursion is stable (if any), it is important to study
the possible existence of minimal solutions. A simple tool that can be used is the Perron-Kreuser
theorem \cite{Kreuser:1914:UTV} (see for instance \cite{Cash:1980:ANO} for a more recent
account of this important result), which in some cases gives simple answers to this question. This is not the case here, however.

To see this, first we notice that the coefficients of the recurrence behave as $k\rightarrow +\infty$
as follows:
\begin{equation}\label{eq:convrep08}
q_1 \sim 16 k^5,\, q_0\sim -16 (4a-b) k^3,\,q_{-1}\sim -8 k^3,\,q_{-2}\sim -4bk ,\, q_{-3}\sim k.
\end{equation}
The coefficients $p_i$ have similar behaviour, and we have $p_i\sim -4kq_i$ as $k\rightarrow +\infty$. Considering
the Newton-Pusieux diagram corresponding to this asymptotic behaviour, we conclude that all the solutions $\beta_k$ of the 
recurrence relation satisfy:
\begin{equation}\label{eq:convrep09}
\limsup_{k\rightarrow +\infty}(k! |\beta_k|)^{1/k}=\frac12,
\end{equation}
where the constant $\frac12$ is the absolute value of the (degenerate) solutions of the characteristic equation
$16\lambda^4-8\lambda^2+1=0$, and the coefficients in this equation come from the coefficients $q_1$, $q_{-1}$ and 
$q_{-3}$. The same result is true for $\alpha_k$. This is all the information the Perron-Kreuser theorem gives in this case;
it is not possible to infer the existence of minimal solutions only with this result. 

\begin{figure}
\begin{center}
\epsfxsize=15cm 
\hspace*{-1.5cm}
\epsfbox{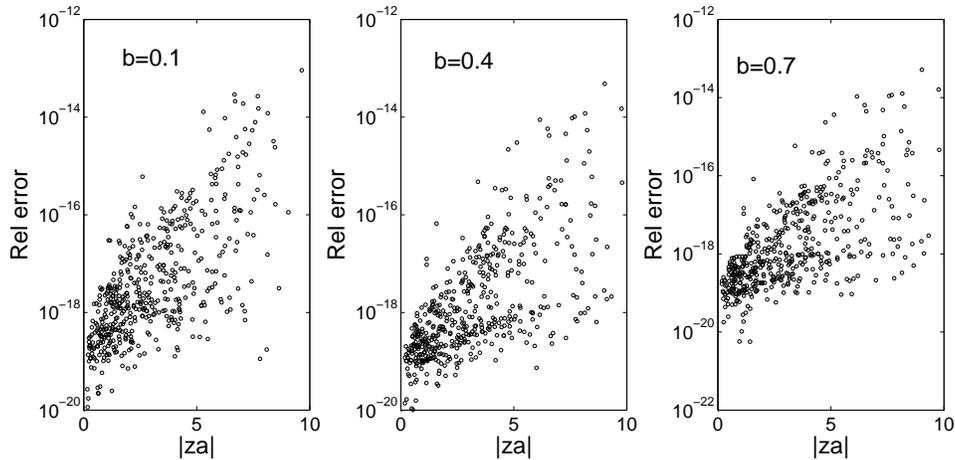}
\caption{Relative errors in the computation of $U(a,b,z)$ using 
the expansions of (\ref{eq:convrep01})
with a fixed number ($N=20$) of $\alpha_k$ and $\beta_k$ coefficients for computing $A_a(z,b)$, $B_a(z,b)$ in (\ref{eq:asrep04}).
\label{fig:fig01}}
\end{center}
\end{figure}

A scaling of the sequences $\alpha_k$ and $\beta_k$ in order to minimize the risk of
overflow/underflow in the recurrence can be used, by defining $\tilde{\alpha}_k=k! 2^k \alpha_k$ and similarly 
for the $\beta_k$'s.
Numerical experiments show that for small $a$, the forward computation appears to be stable while for large $a$
this is no longer true. 

Although, the Perron-Kreuser theorem does not give information on the existence of minimal solutions in this case, it is 
possible to test numerically the existence of minimal solution. One can start the recurrence with large values of $k$
and arbitrarily chosen starting values, apply the recurrence in the backward direction and check whether the same value
$\alpha_1/\alpha_0$ (or $\beta_1 /\beta_0$) is obtained irrespectively of the initial values. If this is the case, this
means that a minimal solution must exist. We have checked  this fact for several sets of starting values and values of
the parameters and it seems that such a minimal solution exists; however, it does not appear to be the requested solution
(because the values do not correspond to (\ref{eq:convrep04})). 

\begin{figure}
\epsfxsize=13cm \epsfbox{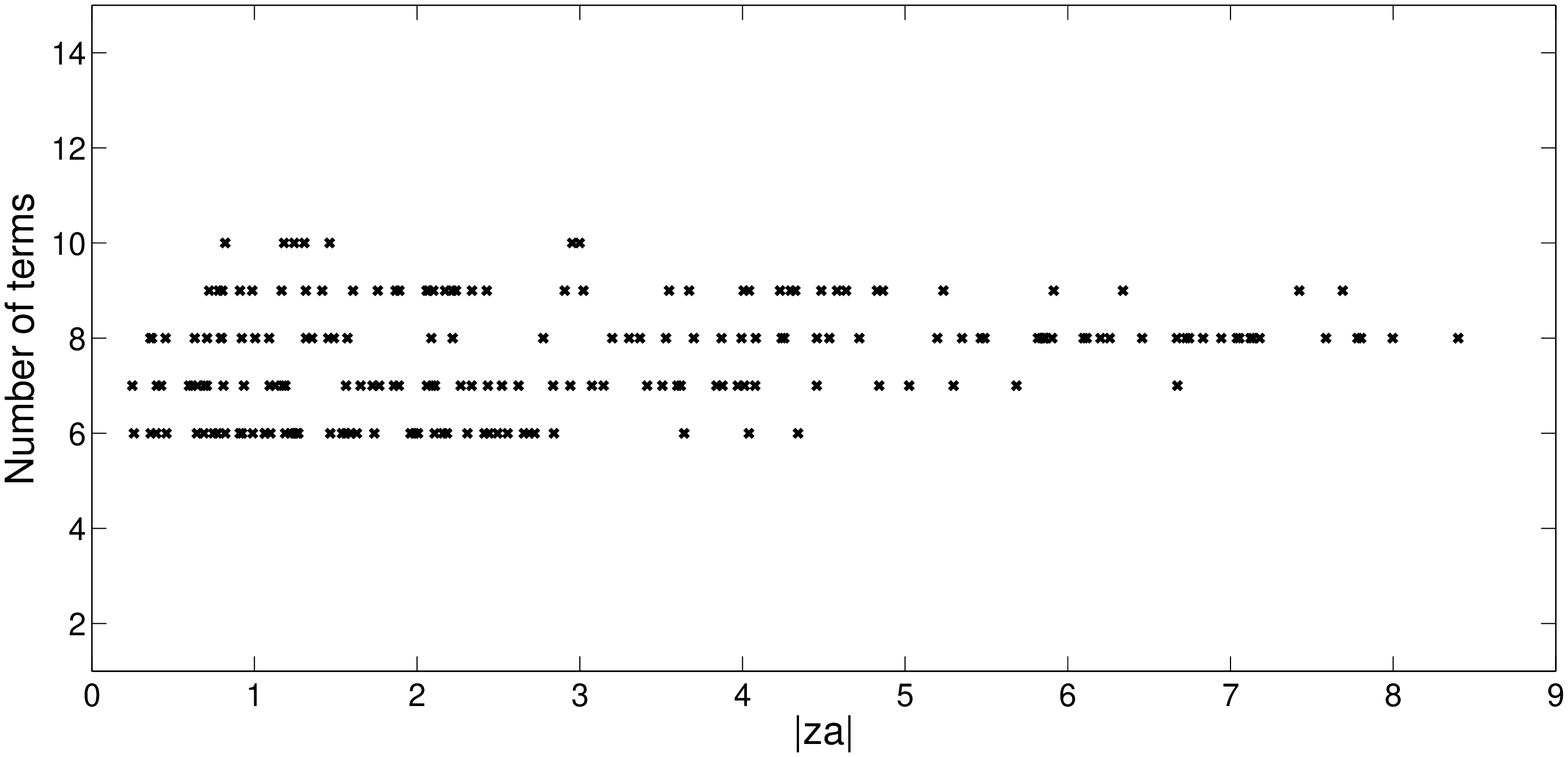}
\caption{Number of terms needed in the expansions of (\ref{eq:convrep01}) to compute $U(a,b,z)$ with an accuracy better than $10^{-14}$.
The value of the parameter $b$ has been fixed to $0.4$.
\label{fig:fig02}}
\end{figure}

Therefore, backward recurrence does not seem to be of help. Contrarily, there is a range of variables for which the
the forward recurrence yields correct results for the expansions. For instance, when $\vert z\vert<1$,  $\vert b\vert<1$ (where the $b$-values 
should be neither too close to 0 nor to 1) we obtain an accuracy close to double precision (or better)
for values of the product $\vert za \vert$ smaller than $10$ for computing the $U$ function, as illustrated in Figure~\ref{fig:fig01}.
In this figure we show, for three values of the $b$ parameter, the relative errors in the computation of $U(a,b,z)$ using 
the expansions of (\ref{eq:convrep01})
with a fixed number  of terms ($N=20$) for computing $A_a(z,b)$, $B_a(z,b)$ in (\ref{eq:convrep01}). On the other hand, when the demanded accuracy 
 in the computation of $U(a,b,z)$ is fixed to $10^{-14}$, for example, the number of terms needed  in (\ref{eq:convrep01}) seems
not to be larger than 10, as shown in Figure~\ref{fig:fig02} for $b=0.4$; this also applies for other values of the
$b$ parameter, excluding values very close to $0$ or $1$.

\section{Acknowledgements}
The authors thank the anonymous referee for helpful comments and suggestions.
The authors acknowledge financial support from 
{\emph{Ministerio de Ciencia e Innovaci\'on}}, project MTM2012-34787. NMT thanks CWI, Amsterdam, for scientific support.

\bibliographystyle{plain}

\def\cprime{$'$}

\end{document}